\documentclass[conference]{IEEEtran}
\IEEEoverridecommandlockouts
\usepackage{cite}
\usepackage{amsmath,amssymb,amsfonts}
\usepackage{algorithmic}
\usepackage{graphicx}
\usepackage{textcomp}
\usepackage{xcolor}
\def\BibTeX{{\rm B\kern-.05em{\sc i\kern-.025em b}\kern-.08em
    T\kern-.1667em\lower.7ex\hbox{E}\kern-.125emX}}

\usepackage{amsthm}
\usepackage{enumerate}
\usepackage{color}
\usepackage[mathscr]{euscript}
\usepackage{url}
\usepackage[all]{xy}
\usepackage{booktabs}
\usepackage{epstopdf}
\usepackage[caption=false]{subfig}
\usepackage{mathtools}

\usepackage{algorithm}
\usepackage{algorithmic}

\usepackage[capitalize,nameinlink,noabbrev]{cleveref}

\usepackage[scale=0.75]{geometry}
\usepackage{tikz}
\usetikzlibrary{shapes}
\usetikzlibrary{arrows}
\usetikzlibrary{matrix}

\newtheorem{theorem}{Theorem}[section]

\newtheorem{lemma}[theorem]{Lemma}

\newtheorem{proposition}[theorem]{Proposition}

\newtheorem{example}[theorem]{Example}

\newcommand{\tens}[1]{\boldsymbol{\mathcal{#1}}}
\newcommand{\tenselem}[1]{\mathcal{#1}}

\newcommand{\matr}[1]{\boldsymbol{#1}}

\newcommand{\vect}[1]{\boldsymbol{#1}}

\newcommand{\set}[1]{\mathscr{#1}}
\newcommand{\T}{{\sf T}}        
\renewcommand{\H}{{\sf H}}      

\makeatletter
\newcommand{\rank}[1]{\mathop{\operator@font rank}\{#1\}}
\newcommand{\colrank}[1]{\mathop{\operator@font colrank}\{#1\}}
\newcommand{\krank}[1]{\mathop{\operator@font krank}\{#1\}}
\newcommand{\trace}[1]{\mathop{\operator@font tr}\{#1\}}
\newcommand{\Diag}[1]{\mathop{\operator@font Diag}\{#1\}}    
\newcommand{\diag}[1]{\mathop{\operator@font diag}\{#1\}}    
\newcommand{\Span}[1]{\mathop{\operator@font Span}\{#1\}}    
\newcommand{\argmin}{\mathop{\operator@font argmin}}

\newcommand{\offdiag}[1]{\mathop{\operator@font offdiag}\{#1\}}    
\newcommand{\Proj}[2]{\mathop{\operator@font Proj_{#1}}{#2}}
\newcommand{\ProjGrad}[2]{\mathop{{\operator@font grad} }#1(#2)}
\newcommand{\Hess}[2]{\mathrm{Hess}_{#2}{#1}}

\makeatother

\newcommand{\eqdef}{\stackrel{\sf def}{=}}
\newcommand{\RR}{\mathbb{R}}
\newcommand{\CC}{\mathbb{C}}

\newcommand{\FF}{\mathbb{F}}

\newcommand{\ON}[1]{\set{O}_{#1}}

\newcommand{\UN}[1]{\set{U}_{#1}}

\newcommand{\Gmat}[3]{\matr{G}^{(#1,#2,#3)}}

\newcommand{\contr}[1]{\mathop{\bullet_{#1}}}   

\newcommand{\ui}{i}

\newcommand{\HessPS}[3]{\mathfrak{D}^{(#1,#2)}_{#3}}

\newcommand{\Utwo}{\Psi}
\newcommand{\hij}[3]{h_{(#1,#2), #3}}
\newcommand{\Gamij}[3]{\matr{\Gamma}^{(#1,#2,#3)}}

\begin{document}

\title{On the convergence of Jacobi-type algorithms for Independent Component Analysis\\
\thanks{This work was supported in part by the National Natural Science Foundation of China (No.11601371), and by  the  Agence Nationale de Recherche (ANR grant LeaFleT, ANR-19-CE23-0021).}
}

\author{\IEEEauthorblockN{Jianze Li}
\IEEEauthorblockA{
\textit{Shenzhen Research Institute of Big Data}\\
Shenzhen, China \\
lijianze@gmail.com}
\and
\IEEEauthorblockN{Konstantin Usevich}
\IEEEauthorblockA{
\textit{CRAN, Universit\'{e} de Lorraine}\\
CNRS, Nancy, France \\
konstantin.usevich@univ-lorraine.fr}
\and
\IEEEauthorblockN{Pierre Comon}
\IEEEauthorblockA{
\textit{GIPSA-Lab, Univ. Grenoble Alpes}\\
CNRS, Grenoble, France \\
pierre.comon@gipsa-lab.fr}
}

\maketitle

\begin{abstract}
Jacobi-type algorithms for simultaneous approximate diagonalization of real (or complex) symmetric tensors have been widely used in independent component analysis (ICA) because of their good performance.
One natural way of choosing the index pairs in Jacobi-type algorithms is the classical cyclic ordering, while the other way is based on the Riemannian gradient in each iteration. 
In this paper, we mainly review in an accessible manner our recent results in a series of papers about weak and global convergence of these Jacobi-type algorithms.
These results are mainly based on the Lojasiewicz gradient inequality. 
\end{abstract}


\section{Introduction}

Let $\mathcal{M}$ be the orthogonal group $\ON{n}$ or the unitary group $\UN{n}$. 
In this paper, we mainly study the following problem, which is to maximize a smooth function 
\begin{equation}\label{eq:min_M}
\textit{f}: \mathcal{M} \to \RR.
\end{equation}
The cost functions in \cref{ex-tensor-diag} and \cref{exam-cost-complex} are all in the form of problem \eqref{eq:min_M},
and they are important in \emph{blind source separation} and \emph{Independent Component Analysis} (ICA) \cite{Como92:elsevier,Como94:sp,de1997signal,Como10:book}.

Let $\FF=\RR$ or $\CC$. 
Let $\FF^{n_1\times\cdots\times n_d}\eqdef\FF^{n_1}\otimes\cdots\otimes\FF^{n_d}$
be the linear space of $d$-th order tensors \cite{kolda2009tensor,comon2014tensors,Cichocki15:review,sidiropoulos2017tensor,qi2017tensor}. 
For $\tens{W}\in\FF^{n\times\cdots\times n}$, 
we denote by $\diag{\tens{W}}$ the vector of its diagonal elements, and $\trace{\tens{W}}$ the trace function. 
We denote by $\|\cdot\|$ the Frobenius norm of a tensor or a matrix,
or the Euclidean norm of a vector.
Operator $\contr{p}$ denotes contraction on the $p$-th index of a tensor; when contracted with a matrix, it is understood that summation is always performed on the second index of the matrix. For instance, $(\tens{A}\contr{1}\matr{M})_{ijk}=\sum_\ell \tenselem{A}_{\ell jk} M_{i\ell}$.

\begin{example}[{\cite[Section 1]{LUC2018}}]\label{ex-tensor-diag}
Let $\{\tens{A}^{(\ell)}\}_{1\leq\ell\leq L}\subseteq\RR^{n\times\cdots\times n}$ be a set of $d$-th order \emph{real symmetric tensors}\footnote{The entries do not change under any permutation of indices.} \cite{Comon08:symmetric}.
The \emph{simultaneous approximate diagonalization of symmetric tensors} problem is to maximize 
\begin{equation}\label{eq:sym_tensor_diagonalization}
f(\matr{Q})= \sum_{\ell=1}^{L} \|\diag{\tens{A}^{(\ell)} \contr{1} \matr{Q}^{\T} \cdots \contr{d} \matr{Q}^{\T}}\|^2,
\end{equation}
where $\matr{Q}\in\ON{n}$. 
This problem includes the following well-known problems in ICA as special cases:\\
(i) approximate tensor diagonalization problem \cite{Como94:ifac,Como94:sp,Lathauwer96:simultaneous,Lathauwer01:ICA}, if $L=1$ and $d>2$;\\
(ii) simultaneous approximate matrix diagonalization problem \cite{Cardoso93:JADE}, if $L>1$ and $d=2$.
\end{example}

\begin{example}[{\cite[Section 2]{ULC2019}}]\label{exam-cost-complex}
Let $\tens{A}^{(\ell)}\in\CC^{n\times\cdots\times n}$ be a $d_\ell$-th order complex tensor for $1\leq\ell\leq L$ (orders $d_1,\ldots, d_L$ are potentially different).
For integers $t_\ell$ satisfying $0 \le t_\ell \le d_\ell$ and  $\alpha_\ell \in \RR$ (possibly negative), the cost function is defined as 
\begin{align}\label{cost-fn-general}
f(\matr{U}) &=\sum\limits_{\ell=1}^{L}
\alpha_\ell \|\diag{\tens{W}^{(\ell)}}\|^2,\ \ \textrm{where}\\
\tens{W}^{(\ell)}&=\tens{A}^{(\ell)} \contr{1} \matr{U}^{\H}\cdots \contr{t_\ell} \matr{U}^{\H} \contr{t_\ell+1} \matr{U}^{\T} \cdots \contr{d_\ell} \matr{U}^{\T}.\notag
\end{align}
Let $d=\max(d_1,\ldots,d_L)$.
It can be shown\footnote{a proof can be found  in  \cite[Section 4]{ULC2019} (see also \cite[Prop. 3.5]{jiang2016characterizing}).} that $f$ admits representation \eqref{cost-fn-general} if and only if there exists a $2d$-th order \emph{Hermitian}\footnote{it means that $\tenselem{B}_{i_1\cdots i_d j_1 \cdots j_d} = \tenselem{B}^{*}_{ j_1 \cdots j_d i_1\cdots i_d}.$}\cite{nie2019hermitian} complex tensor $\tens{B}$, 
such that $f$ has a representation
\begin{align}\label{cost-fn-general-trace}
 f(\matr{U}) &= \trace{\tens{V}},\ \ \textrm{where}\\
\tens{V}&= \tens{B} \contr{1} \matr{U}^{\H} \cdots \contr{d} \matr{U}^{\H}\contr{d+1} \matr{U}^{\T} \cdots \contr{2d} \matr{U}^{\T}. \notag
\end{align}
This cost function includes the following well-known problems in ICA as special cases: \\
(i) \emph{simultaneous approximate Hermitian diagonalization of complex matrices} \cite{cardoso1996jacobi}. Let $\matr{A}^{(\ell)} \in \CC^{n\times n}$. 
The cost function is defined as 
\begin{equation}\label{eq-cost-jade}
f(\matr{U}) =  \sum\limits_{\ell=1}^{L}\|\diag{\matr{U}^{\H}\matr{A}^{(\ell)}\matr{U}}\|^2.
\end{equation}
(ii) \emph{approximate diagonalization of a 3rd order complex tensor} \cite{de1997signal,Como10:book}.
Let $\tens{A}\in\CC^{n\times n\times n}$.
The cost function is defined as 
\begin{equation}\label{eq-cost-3-tensor}
f(\matr{U}) =  \|\diag{\tens{A} \contr{1} \matr{U}^{\H} \contr{2} \matr{U}^{\T}\contr{3} \matr{U}^{\T}}\|^2. 
\end{equation}
(iii) \emph{approximate  diagonalization of a 4th order complex tensor} \cite{Como04:ijasp}. 
Let $\tens{B}\in\CC^{n\times n\times n\times n}$ be 
Hermitian. 
The cost function is defined as 
\begin{equation}\label{cost-4-order}
f(\matr{U}) = \trace{\tens{B} \contr{1} \matr{U}^{\H} \contr{2} \matr{U}^{\H}\contr{3} \matr{U}^{\T}\contr{4} \matr{U}^{\T}}. 
\end{equation}
\end{example}

To solve problem \eqref{eq:min_M},
Jacobi-like sweeping procedure is a popular method,
which can break a high-dimensional optimization problem into a sequence of one or two dimensional subproblems.
This method often has a very high speed,
as the solution of the subproblem can often be written in closed form.
In particular, there have been several Jacobi-type algorithms proposed for the cost functions \eqref{eq:sym_tensor_diagonalization} and
\eqref{cost-fn-general} in ICA,
\emph{e.g.}, the
well-known Jacobi CoM2 algorithm \cite{Como92:elsevier,Como94:sp,Como94:ifac,Como10:book}
and \emph{Joint approximate diagonalization of eigenmatrices} (JADE) algorithm \cite{Cardoso93:JADE,cardoso1996jacobi}.

In this paper, we mainly review our recent results   \cite{LUC2018,LUC2018-LAA,ULC2019,LUC2019} about the \emph{weak convergence}\footnote{every accumulation point is a stationary point.} and \emph{global convergence}\footnote{for any starting point, the iterates always converge to a limit point.} of Jacobi-type algorithms on $\mathcal{M}$.
In \cref{sec-jaco-alg}, we introduce the general Jacobi algorithm on $\mathcal{M}$, as well as the Jacobi-C algorithms.
In \cref{sec:Jacobi_G}, we introduce the Jacobi-G algorithms, and present our convergence results about them.
\cref{sec:experennts} includes some experiments. 

\section{Jacobi-type algorithms}\label{sec-jaco-alg}

\subsection{General Jacobi algorithm on $\ON{n}$}
Let $\theta\in\RR$ be an angle and $(i,j)$ be a pair of indices with $1 \leq i < j \leq n$.
Denote by $\Gmat{i}{j}{\theta}$ the \emph{Givens rotation} matrix, as defined \textit{e.g.}, in \cite{GoluV96:jhu,Como10:book,LUC2018}.
The \emph{general Jacobi algorithm} for problem \eqref{eq:min_M} on $\ON{n}$ can be summarized as in \cref{alg:jacobi}.
In this algorithm, if several equivalent maximizers are present in \eqref{eq:h_k}, we choose the one with the angle of smallest magnitude.

\begin{algorithm}\caption{General Jacobi algorithm on $\ON{n}$}\label{alg:jacobi}
\begin{algorithmic}[1]
	\STATE {\bf Input:} a starting point $\matr{Q}_{0}$
	\STATE {\bf Output:} sequence of iterates $\matr{Q}_{k}$
	\FOR{ $k=1,2,\ldots$}
	\STATE Choose the pair $(i_k,j_k)$ according to a certain pair selection rule
	\STATE Compute $\theta^{*}_{k}$ that maximizes
	\begin{equation}\label{eq:h_k}
	h_k(\theta)\eqdef{f}(\matr{Q}_{k-1} \Gmat{i_k}{j_k}{\theta})
	\end{equation}
	\STATE Update $\matr{Q}_k = \matr{Q}_{k-1}\Gmat{i_k}{j_k}{\theta^{*}_k}$ 
	\ENDFOR
\end{algorithmic}
\end{algorithm}

\subsection{General Jacobi algorithm on $\UN{n}$}

Let $\matr{\Utwo} \in \UN{2}$. 
Denote by $\Gmat{i}{j}{\matr{\Utwo}}$ the \emph{plane transformation} matrix, as defined \textit{e.g.}, in \cite{ULC2019,Como10:book}.
The \emph{general Jacobi algorithm} for problem \eqref{eq:min_M} on $\UN{n}$ can be summarized as in \cref{Jacobi-general}.

\begin{algorithm}\caption{General Jacobi algorithm on $\UN{n}$}\label{Jacobi-general}
\begin{algorithmic}[1]
	\STATE {\bf Input:} a starting point $\matr{U}_{0}$
	\STATE {\bf Output:} sequence of iterates $\matr{U}_{k}$
	\FOR{$k=1,2,\ldots$}
	\STATE Choose the pair $(i_k,j_k)$ according to a certain pair selection rule
	\STATE Compute $\matr{\Utwo}_k^{*}$ that maximizes \begin{equation}\label{eq-func-h-new}
	h_k(\matr{\Utwo})\eqdef f(\matr{U}_{k-1}\Gmat{i_k}{j_k}{\matr{\Utwo}})
	\end{equation}
	\STATE Update $\matr{U}_k = \matr{U}_{k-1} \Gmat{i_k}{j_k}{\matr{\Utwo}_k^{*}}$
	\ENDFOR
\end{algorithmic}
\end{algorithm}


\subsection{Jacobi-C algorithms}

One natural pair selection rule in \cref{alg:jacobi} and \cref{Jacobi-general} is in cyclic fashion \cite{GoluV96:jhu,Como10:book} as follows: 
\begin{equation}\label{equation-Jacobi-C}
{\small \begin{split}
&(1,2) \to (1,3) \to \cdots \to (1,n) \to \\
& (2,3) \to \cdots \to (2,n) \to \\
& \cdots  \to (n-1,n) \to \\
&(1,2) \to (1,3) \to \cdots.
\end{split}}
\end{equation}

The Jacobi algorithm with cyclic rule \eqref{equation-Jacobi-C} is called \emph{Jacobi-C} algorithm.
Although this cyclic rule is very simple, the convergence of Jacobi-C algorithms is difficult to study.
In \cite[Remark 6.5]{LUC2018-LAA}, 
we proved that, if $f$ is the cost function in \eqref{eq:sym_tensor_diagonalization} with $d=3$ and $L=1$, then Jacobi-C algorithm may converge to a saddle point of $f$.
To our knowledge,
the convergence of Jacobi-C algorithms for cost functions \eqref{eq:sym_tensor_diagonalization} and
\eqref{cost-fn-general} is still unknown, except in the single matrix case \cite{GoluV96:jhu}. 

\section{Convergence analysis of Jacobi-G algorithms}\label{sec:Jacobi_G}

In this section,
we mainly present our recent results  \cite{LUC2018,ULC2019,LUC2019} about weak and global convergence of \cref{alg:jacobi} and \cref{Jacobi-general} under a gradient based pair selection rule. 
Before that, we need to first recall the well-known \emph{\L{}ojasiewicz gradient inequality}\cite{law1965ensembles,AbsMA05:sjo,SU15:pro,bolte2014proximal}.

\subsection{\L{}ojasiewicz gradient inequality}

Let ${T}_{\vect{x}} \mathcal{M}$ be the tangent space at $\vect{x}$. 
Let $\nabla\textit{f}(\vect{x})$ be the Euclidean gradient, and $\ProjGrad{f}{\vect{x}}$ be the  projection of $\nabla\textit{f}(\vect{x})$ on ${T}_{\vect{x}} \mathcal{M}$, \emph{i.e.}, the Riemannian gradient \cite{Absil08:Optimization}. 
The following results were proved in \cite{SU15:pro}.

\begin{lemma}\label{lemma-SU15}
Let $\mathcal{M}\subseteq\RR^n$ be an analytic\footnote{See {\cite[Definition 2.7.1]{krantz2002primer}} or \cite[Definition 5.1]{LUC2018} for a definition of an analytic submanifold.} submanifold.
Then any point $\vect{x}\in \mathcal{M}$ satisfies a \L{}ojasiewicz inequality for
$\ProjGrad{f}{\cdot}$, 
that is,
there exist $\delta>0$, $\sigma>0$ and $\zeta\in (0,1/2]$
such that for all $\vect{y}\in\mathcal{M}$ with $\|\vect{y}-\vect{x}\|<\delta$,
it holds that
\begin{equation}\label{eq-kl-inequa}
|\textit{f}(\vect{y})-\textit{f}(\vect{x})|^{1-\zeta}\leq \sigma\|\ProjGrad{f}{\vect{y}}\|.
\end{equation}
\end{lemma}

\begin{theorem}[{\cite[Theorem 2.3]{SU15:pro}}]\label{theorem-SU15}
Let $\mathcal{M}\subseteq\RR^n$ be an analytic submanifold
and $\{\vect{x}_k\}_{k\geq 1}\subseteq\mathcal{M}$ be a sequence.
Suppose that $f$ is real analytic and, for large enough $k$,\\
(i) there exists $\sigma>0$ such that
$$|\textit{f}(\vect{x}_{k+1})-\textit{f}(\vect{x}_k)|\geq \sigma\|\ProjGrad{f}{\vect{x}_k}\|\|\vect{x}_{k+1}-\vect{x}_{k}\|;$$
(ii) $\ProjGrad{f}{\vect{x}_k}=0$ implies that $\vect{x}_{k+1}=\vect{x}_{k}$.\\
Then any accumulation point $\vect{x}_*$ of $\{\vect{x}_k\}_{k\geq 1}\subseteq\mathcal{M}$ is the only limit point.

If, in addition, for some $\kappa > 0$  and for large enough $k$ it holds that 
\begin{equation*}\label{eq:safeguard}
\|\vect{x}_{k+1} - \vect{x}_{k}\| \ge \kappa \|\ProjGrad{f}{\vect{x}_k}\|,
\end{equation*}
then the following convergence rates apply
\begin{equation*}
\|\vect{x}_{k} - \vect{x}_*\| \le  C\begin{cases}
e^{-ck}, & \text{ if } \zeta = \frac{1}{2} \text{ (for some }c > 0\text{)}, \\
k^{-\frac{\zeta}{1-2\zeta}}, & \text{ if } 0 < \zeta < \frac{1}{2},
\end{cases}
\end{equation*}
where $\zeta$ is the parameter in \eqref{eq-kl-inequa} at the limit point $\vect{x}_*$.
\end{theorem}

\subsection{Convergence of Jacobi-G algorithm on $\ON{n}$}

A gradient based pair selection rule of \cref{alg:jacobi} was proposed in \cite{IshtAV13:simax},  which chooses  a pair $(i_k,j_k)$ in each iteration satisfying that 
\begin{equation}\label{eq:pair_selection_gradient}
|h'_k(0)| \ge \delta \|\ProjGrad{f}{\matr{Q}_{k-1}}\|,
\end{equation}
where  $\delta$ is a small positive constant.
The main idea behind this rule is to choose Givens rotations that are well aligned with the Riemannian gradient of $f$. 
We call the Jacobi algorithm with this rule the \emph{Jacobi-G} algorithm, which is summarized in \cref{alg:jacobi-G}. 

\begin{algorithm}\caption{Jacobi-G algorithm on $\ON{n}$}\label{alg:jacobi-G}
\begin{algorithmic}[1]
	\STATE {\bf Input:} a starting point $\matr{Q}_{0}$, a positive constant $0<\delta<{2}/n$
	\STATE {\bf Output:} sequence of iterates $\matr{Q}_{k}$
	\FOR{ $k=1,2,\ldots$}
	\STATE Choose the pair $(i_k,j_k)$ satisfying \eqref{eq:pair_selection_gradient}
	\STATE Compute $\theta^{*}_{k}$ that maximizes the function \eqref{eq:h_k}
	\STATE Update $\matr{Q}_k = \matr{Q}_{k-1}\Gmat{i_k}{j_k}{\theta^{*}_k}$ 
	\ENDFOR
\end{algorithmic}
\end{algorithm}

It was shown in \cite[Lemma 5.2]{IshtAV13:simax} and \cite[Lemma 3.1]{LUC2018} that, in \cref{alg:jacobi-G}, we can always choose such a pair $(i_k,j_k)$ satisfying condition \eqref{eq:pair_selection_gradient}.
Then, based on the proof of \cite[Theorem 5.4]{IshtAV13:simax}, we can easily get the following result about the weak convergence of \cref{alg:jacobi-G}. 

\begin{theorem}[{\cite[Theorem 3.3]{LUC2018}}]\label{theorem-ishteva-stationary-point}
Let $f$ be smooth. Then every accumulation point of the iterates in \cref{alg:jacobi-G} is a stationary point of $f$.
\end{theorem}

Based on \cref{theorem-SU15}, the following global convergence result of \cref{alg:jacobi-G} was proved.

\begin{theorem}[{\cite[Theorem 5.6]{LUC2018}}]\label{theorem-convergence-order-3}
Let $f$ be the cost function in \eqref{eq:sym_tensor_diagonalization} with $d=2$ or $d=3$.
Then, for any starting point $\matr{Q}_{0}$, the iterates in \cref{alg:jacobi-G} always converge to a  stationary point of $f$.
\end{theorem}

\subsection{Convergence of Jacobi-G algorithm on $\UN{n}$}

As a complex generalization of \cref{alg:jacobi-G},
the following \emph{Jacobi-G algorithm on $\UN{n}$} was formulated in \cite[Section 2]{ULC2019}. 

\begin{algorithm}\caption{Jacobi-G algorithm on $\UN{n}$}\label{Jacobi-G-general}
\begin{algorithmic}[1]
	\STATE {\bf Input:} a starting point $\matr{U}_{0}$, a positive constant $0<\delta<\sqrt{2}/n$
	\STATE {\bf Output:} sequence of iterates $\matr{U}_{k}$
	\FOR{$k=1,2,\ldots$}
	\STATE Choose an index pair $(i_k,j_k)$ satisfying
	\begin{equation}\label{inequality-gra-based}
	\|\ProjGrad{h_k}{\matr{I}_2}\| \geq \delta
	\|\ProjGrad{f}{\matr{U}_{k-1}}\|
	\end{equation}
	\STATE Compute $\matr{\Utwo}^{*}_k$ that maximizes the function \eqref{eq-func-h-new}
	\STATE Update $\matr{U}_k = \matr{U}_{k-1} \Gmat{i_k}{j_k}{\matr{\Utwo}_k^{*}}$
	\ENDFOR
\end{algorithmic}
\end{algorithm}

It was shown in \cite[Section 4]{ULC2019} that, in \cref{Jacobi-G-general}, we can always choose such a pair $(i_k,j_k)$ satisfying condition \eqref{inequality-gra-based}.
Then, based on \cite[{Theorem 2.5}]{boumal2016global},
the following weak convergence result of \cref{Jacobi-G-general} was proved in \cite[Section 5]{ULC2019}. 

\begin{proposition}\label{theorem-local-gene}
Let $f$ be the cost function in \eqref{cost-fn-general}. 
Then every accumulation point of the iterates in \cref{Jacobi-G-general} is a stationary point of $f$. 
\end{proposition}

For $\matr{U}\in\UN{n}$ and $1 \le i < j \le n$, we define the restriction of $f$ as 
\begin{equation*}
\hij{i}{j}{\matr{U}}: \UN{2} \longrightarrow \RR,\ 
{\matr{\Utwo}}\longmapsto f(\matr{U}\Gmat{i}{j}{\matr{\Utwo}}).
\end{equation*}
When maximizing $\hij{i}{j}{\matr{U}}$, we only\footnote{this is reasonable; see \cite[Section 2]{ULC2019} for more discussions.} consider 
\begin{equation*}
\matr{\Psi}= \matr{\Utwo}(c,s_1,s_2)
= \begin{bmatrix}
c & -(s_1+\ui s_2) \\
s_1-\ui s_2 & c
\end{bmatrix},
\label{unitary-para-2}
\end{equation*}
satisfying $c^2+s_1^2+s_2^2 =1$. 
It was proved in \cite[Section 4]{ULC2019} that, if $f$ is the cost function in \eqref{cost-fn-general-trace} with $d\leq3$, 
there exists a real symmetric matrix $\Gamij{i}{j}{\matr{U}}\in\RR^{3\times 3}$ such that
\begin{equation}\label{eq-cost-quadratic-form}
\hij{i}{j}{\matr{U}}(c,s_1,s_2) = \vect{r}^{\T} \Gamij{i}{j}{\matr{U}}\vect{r},
\end{equation}
where $\vect{r} \eqdef (2c^2-1, -2cs_1, -2cs_2)^{\T}.$
Therefore, the cost functions \eqref{eq-cost-jade},  \eqref{eq-cost-3-tensor} and \eqref{cost-4-order} in \cref{exam-cost-complex} all satisfy \eqref{eq-cost-quadratic-form}. 
Now, we define 
\begin{equation*}
\HessPS{i}{j}{\matr{U}}  \eqdef  2\left(\begin{bmatrix}\Gamij{i}{j}{\matr{U}}_{2,2} & \Gamij{i}{j}{\matr{U}}_{2,3} \\ \Gamij{i}{j}{\matr{U}}_{3,2} & \Gamij{i}{j}{\matr{U}}_{3,3}  \end{bmatrix} - \Gamij{i}{j}{\matr{U}}_{1,1}\matr{I}_2\right).
\end{equation*}

Let $\Hess{f}{\vect{x}}$ be the Riemannian Hessian\footnote{see \cite[p.105]{Absil08:Optimization} for a detailed definition.} of $f$ at $\vect{x}\in\mathcal{M}$, which is a linear map $\mathbf{T}_{\vect{x}}\mathcal{M} \to \mathbf{T}_{\vect{x}}\mathcal{M}$. 
Based on \cref{theorem-SU15}, the following global convergence result of \cref{Jacobi-G-general} was proved in \cite[Section 7]{ULC2019}.

\begin{theorem}\label{thm:accumulation_points}
Let $f$ be the cost function in \eqref{cost-fn-general-trace} with $d \le 3$. 
Let $\matr{U}_{*}$ be  an accumulation point of the iterates in  \Cref{Jacobi-G-general} (and thus $\ProjGrad{f}{\matr{U}_{*}} = 0$ by \Cref{theorem-local-gene}). 
Assume that  $\HessPS{i}{j}{\matr{U}_{*}}$ is  negative definite for all pairs $(i,j)$.
Then\\
(i) $\matr{U}_{*}$ is the only limit  point and convergence rates in \Cref{theorem-SU15} apply.\\
(ii)  If the rank of Riemannian Hessian is maximal at $\matr{U}_{*}$ (\emph{i.e.}, $\rank{\Hess{f}{\matr{U}_{*}}} =  n(n-1)$), then the speed of convergence is linear.
\end{theorem}

\section{Experiments}\label{sec:experennts}

In this section, we make some experiments to see the convergence behaviours of Jacobi-C and Jacobi-G algorithms on $\ON{n}$. 

\begin{example}\rm\label{example-2}
We randomly generate two symmetric tensors $\tens{A}_{1}\in\RR^{10\times 10\times 10}$ and 
$\tens{A}_{2}\in\RR^{10\times 10\times 10\times 10}$. 
Let $f$ be the cost function in \eqref{eq:sym_tensor_diagonalization} with $L=1$.  
Let the starting point $\matr{Q}_{0}=\matr{I}_{10}$ in Jacobi-C algorithm and \cref{alg:jacobi-G}, and $\delta=0.1$ in \cref{alg:jacobi-G} . 
The results are shown in \cref{figure-example-1}.

\begin{figure}[tbhp]
\centering
\subfloat[$\tens{A}_1\in\RR^{10\times 10\times 10}$]{\includegraphics[width=0.425\textwidth]{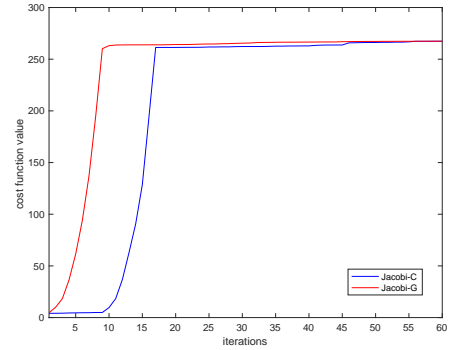}}\!\!\!
\subfloat[$\tens{A}_2\in\RR^{10\times 10\times 10\times 10}$]{\includegraphics[width=0.425\textwidth]{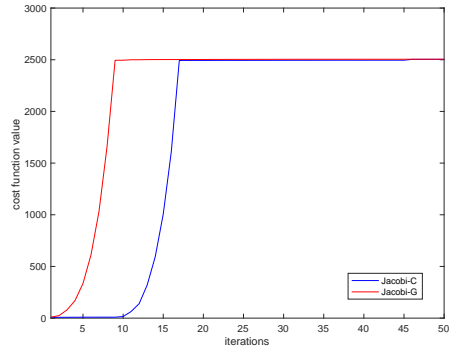}}
\caption{Results of Jacobi-C and Jacobi-G algorithms on $\ON{n}$.}
\label{figure-example-1}
\end{figure}

\end{example}

\section{Conclusion}

In this paper, we mainly review our results  \cite{LUC2018,LUC2018-LAA,ULC2019,LUC2019} about the convergence of Jacobi-type algorithms on $\ON{n}$ and $\UN{n}$.
For the moment, there are still some open problems: \\
(i) If $f$ is the cost function in \eqref{eq:sym_tensor_diagonalization} with $d\geq 4$, the global convergence of \cref{alg:jacobi-G} is still unknown.\\
(ii) The global convergence of Jacobi-C algorithm is still unknown.

\bibliographystyle{IEEEtran}
\bibliography{Ref_SAM2020}

\end{document}